\documentclass[10pt, a4paper]{article}

\usepackage{fancyhdr}
\usepackage{color}
\usepackage{titletoc}
\usepackage{latexsym}
\usepackage{latexsym,bm}
\usepackage{amsmath}
\usepackage{amssymb}
\usepackage{multicol}
\usepackage{graphics}
\usepackage{graphicx}
\usepackage{subfigure}
\usepackage{indentfirst}
\usepackage{epsfig}
\usepackage{mathrsfs}

\begin{document}
\title{{\large  \bf The Extension Degree Conditions for Fractional Factor}}
\author{{\large  Wei Gao$^{1}$, Weifan Wang$^{2}$, Juan L.G. Guirao$^{3}$\thanks{Corresponding author: juan.garcia@upct.es}}
\\
{\small 1. School of Information Science and Technology, Yunnan Normal University, Kunming 650500, China}\\
{\small 2. Department of Mathematics, Zhejiang Normal University,
Jinhua 321004, China}\\
{\small 3. Departamento de Matem\'{a}tica Aplicada y Estad\'{\i}stica,
Universidad Polit\'{e}cnica de Cartagena, }\\
{\small Hospital de Marina,
30203-Cartagena, Regi\'{o}n de Murcia, Spain}\\
\date{}
}

\maketitle
\newtheorem{claim}{Claim}
\newtheorem{problem}{Problem}
\newtheorem{theorem}{Theorem}
\newtheorem{lemma}{Lemma}
\newtheorem{corollary}{Corollary}
\newtheorem{definition}[theorem]{Definition}
\newtheorem{proposition}[theorem]{Proposition}
\newtheorem{conjecture}{Conjecture}
\newtheorem{remark}{Remark}
\newcommand{\no}{\noindent}
\newcommand{\qed}{\hfill $\Box$ }

\parindent=0.5cm

\noindent  {\bf Abstract:}\ \  In Gao's previous work, the
authors determined several graph degree conditions of a graph which admits fractional factor in particular settings. It was revealed
that these degree conditions are tight if $b=f(x)=g(x)=a$ for all vertices
$x$ in $G$. In this paper, we continue to discuss these degree
conditions for admitting fractional factor in the setting that several vertices and edges are removed and there is
a difference $\Delta$ between $g(x)$ and $f(x)$ for every vertex
$x$ in $G$. These obtained new degree conditions  reformulate Gao's previous conclusions, and show how $\Delta$ acts in the results. Furthermore, counterexamples are structured to reveal the sharpness of degree conditions in the setting $f(x)=g(x)+\Delta$.



\section{Introduction}
In many engineering applications, their mathematical models can be
expressed as a (direct or undirect) graph. For example, we look upon the
network as a graph. Some correspondences are given here: the site matches with a vertex and the channel matches with an edge in the graph. In
conventional network, the mission of data transmission is based on the selection of
the shortest way between vertices. However, the
computation of network flow in software
definition network determines the data transmission. It chooses the path that is least congested at present. In this
view, the pattern of data transmission problem  in SDN setting is just the
existence of fractional factor in the corresponding graph.

Graph $G=(V,E)$ mentioned here are all simple graph
with its edge set $E(G)$ and its vertex set $V(G)$. Throughout
this paper, we set $n=|V(G)|$ as the order of graph. For a vertex
$x$ in $G$,
  $N_{G}(x)$ and $d_{G}(x)$ are used to denote  the
neighborhood and the degree of $x$ in $G$, respectively. Let
$N_{G}[x]=N_{G}(x)\cup \{x\}$. To simplify, we use $N(x)$,$d(x)$
and $N[x]$ to express $d_{G}(x)$, $N_{G}(x)$ and $N_{G}[x]$,
respectively. Set
$\delta(G)$ as the minimum degree of $G$. We set $G[S]$ as the sub-graph of
$G$ deduced from $S\subseteq V(G)$, and $G-S=G[V(G) \backslash
S]$. Set $e_{G}(S_{1},S_{2})=|\{e=v_{1}v_{2}|v_{1}\in S_{1},v_{2}\in S_{2}\}|$ for any $S_{1},S_{2}\subseteq V(G)$ with $S_{1}\cap S_{2}=\emptyset$. Denote
$\sigma_{2}(G)=\min\{d_{G}(u)+d_{G}(v)|uv\notin E(G)\}$. The other
terms used without clear definitions here can be refered to classic graph theory book \cite{Bondy08}.

Functions $f$ and $g$ are two integer-valued defined on $V(G)$
satisfying $f(x)\ge g(x)\ge0$ for all vertices $x$ in $G$. A {\em fractional
$(g,f)$-factor} is regarded as a score function $h$ which maps to
every element in $E(G)$ a real number belongs to [0,1]. As a result, for every
vertex $x$ we get $g(x)\le d_{G}^{h}(x)\le f(x)$,  and
$\sum\limits_{e\in E(x)}h(e)=d_{G}^{h}(x)$  where $E(x)=\{y|yx\in E(G)\}$. Fractional $f$-factor is regarded as a special case of fractional $(g,f)$-factor
 if the values of two functions are equal for any vertex $x$ in $G$. Fractional $[a,b]$-factor is another special case of
fractional $(g,f)$-factor if
$f(x)=b$, $g(x)=a$ for any vertex $x$ in $G$. In addition,  if the value of both $f$ and $g$  equal to
$k\in\Bbb N$ for any vertex $x$ in $G$, then it's a fractional $k$-factor.

A fractional $(g,f,m)$-deleted graph and a fractional $(g,f,n')$-critical graph imply the existence of fractional factor in special setting when delete $m$ edges and $n'$ vertices, respectively. As the combination of the above two concepts, Gao \cite{Gao12} introduced fractional $(g,f,n',m)$-critical deleted graph to denote  a graph to be
fractional $(g,f,m)$-deleted after removing
any $n'$ vertices. When functions $g$ and $f$ take special value for all vertices, the fractional $(g,f,n',m)$-critical deleted graph becomes various names which are presented in Table \ref{table1}.
\begin{table}
\caption{Special cases of fractional $(g,f,n',m)$-critical deleted graph}\label{table1}
\begin{tabular}{ll}
\hline
setting (for any $v\in V(G)$) & name \\
\hline
$g(x)=f(x)$ & fractional $(f,n',m)$-critical
deleted graph\\
\hline
$f(x)=b$ and $g(x)=a$ &fractional $(a,b,n',m)$-critical
deleted graph\\
\hline
$f(x)=g(x)=k$ & fractional
$(k,n',m)$-critical deleted graph\\
\hline
\end{tabular}
\end{table}
Several recent contributions
in this topic were presented in
Zhou et al.  \cite{5}, \cite{4}, \cite{7}, \cite{6} and \cite{8}, and Gao et al. \cite{add1}, \cite{add2}, \cite{add4} and \cite{add3}, Knor et al. \cite{AMNS1}, and Liu et al. \cite{AMNS2}.

In Zhou \cite{Zhou15} and Zhou et al. \cite{zhouetalinpress}, the
setting was different from the previous situations in which there
is a difference $\Delta$ between  $g(x)$ and $f(x)$ for every vertex
$x$ in $G$, i.e., $b-\Delta\ge f(x)-\Delta\ge g(x)\ge a$ for every
$x$ in $G$. We observe that if $\Delta=0$, then binding number (minimum $\frac{|N(X)|}{|X|}$ where
$\emptyset\ne X\subset V(G)$) condition for
ID-$(g,f)$-factor-critical graph (this concept will be explain later) is
$$bind(G)>\frac{(n-1)(b+2a-1)}{an-(b+a-2)}.$$
After adding the variable $\Delta$ (i.e., $a\le g(x)\le
f(x)-\Delta\le b-\Delta$), by the conclusion obtained by Zhou et
al. \cite{zhouetalinpress}, the binding number condition becomes
$$bind(G)>\frac{(n-1)(b+2a-1+\Delta)}{(a+\Delta)n-(b+a-2)}.$$
This fact reveals that if the setting changes, the lower bound of
binding number for ID-$(g,f)$-factor-critical graph is changed as
well, and the new binding number heavily depends on $\Delta$. There is one thing we must emphasize here is that all the results in this paper are independent from the maximum degree of the graph, and $\Delta$ is only used to represent the difference between $g$ and $f$ throughout the article.

In our article, we consider the degree condition for the existence of fractional factors in the setting that $b-\Delta\ge f(x)-\Delta\ge g(x)\ge a$ for every vertex $x$.  Intuitively, in
our new setting, the new degree conditions should be relied on the
variable $\Delta$, or at least the new degree conditions are
different from the previous ones. Thus, it inspired us to strictly study it theoretically.

In the following context, we first present the major results of part one in fractional $(g,f,n',m)$-critical deleted setting and prove it in details which extended
Theorem 1-3 raised in Gao et. al.
\cite{Gaowei4}, perspectively.
\begin{theorem}\label{theorem1}
Assume $G$ is a graph with $n$ vertices, and set $b,a,n',m$, and $\Delta$
as non-negative integers meeting $b-\Delta\ge a\ge2$ and $n>
\frac{(b+a+2m-2)(b+a)}{\Delta+a}+n'$. Functions $g,f$ are integer-valued
on its vertex set and $b-\Delta\ge  f(x)-\Delta\ge g(x)\ge a$ for every vertex $x$. Then $G$ is fractional
$(g,f,n',m)$-critical deleted if
$\delta(G)\ge\frac{(b-\Delta)n+(\Delta+a)n'}{b+a}$.
\end{theorem}

\begin{theorem}\label{theorem2}
Assume $G$ is a graph with $n$ vertices, and set $b,a, n',m$, and $\Delta$
as non-negative integers meeting $b-\Delta\ge a\ge2$,
$\delta(G)\ge m+n'+\frac{b(b-\Delta)}{\Delta+a}$ and $n>
\frac{(b+a+2m-1)(a+b)}{\Delta+a}+n'$. Functions $g,f$ as integer-valued on its vertex set and meet $b-\Delta\ge f(x)-\Delta\ge g(x)\ge a$ for every vertex $x$ in $G$.
Then $G$ is fractional $(g,f,n',m)$-critical deleted if for any $xy\ne E(G)$, we have
$$\max\{d_{G}(x),d_{G}(y)\}\ge\frac{(b-\Delta)n+(\Delta+a)n'}{b+a}.$$
\end{theorem}

\begin{theorem}\label{theorem3}
Assume $G$ is a graph with $n$ vertices, and set $b,a, n',m$, and $\Delta$
as non-negative integers meeting $b-\Delta\ge a\ge2$,
$\delta(G)\ge\frac{b(b-\Delta)}{\Delta+a}+m+n'$ and $n>
\frac{(b+a+2m-2)(a+b)}{\Delta+a}+n'$. Functions $g,f$ are
integer-valued defined on the vertex set so that $b-\Delta\ge f(x)-\Delta\ge g(x)\ge a$ for every vertex $x$ in $G$. Then $G$ is fractional
$(g,f,n',m)$-critical deleted  if
$\sigma_{2}(G)\ge\frac{2(n(b-\Delta)+n'(\Delta+a))}{b+a}$.
\end{theorem}

The above three theorems manifest conditions for a graph to be fractional
$(g,f,n',m)$-critical deleted from different
aspects. The corollaries on
fractional $(g,f,m)$-deleted graphs can be stated in Table \ref{table2}.
\begin{table}
\caption{Three degree conditions of fractional $(g,f,m)$-deleted graph by setting  $n'=0$ in above three theorems}\label{table2}
\begin{tabular}{lll}
\hline
order of graph & degree condition& additional condition\\
\hline
$n>
\frac{(b+a+2m-2)(a+b)}{\Delta+a}$ & $\delta(G)\ge\frac{(b-\Delta)n}{b+a}$&  \\
\hline
$n>
\frac{(b+a+2m-1)(b+a)}{\Delta+a}$ &$\max\{d_{G}(x),d_{G}(y)\}\ge\frac{(b-\Delta)n}{b+a}$ & $\delta(G)\ge\frac{b(b-\Delta)}{\Delta+a}+m$\\
\hline
$n>
\frac{(b+a+2m-2)(b+a)}{\Delta+a}$ & $\sigma_{2}(G)\ge\frac{2(b-\Delta)n}{b+a}$& $\delta(G)\ge\frac{b(b-\Delta)}{\Delta+a}+m$\\
\hline
\end{tabular}
\end{table}

The data in Table \ref{table2} can be regarded as the extension of Corollary 1,
Corollary 2 and Corollary 3 in Gao et al. \cite{Gaowei4},
respectively. Furthermore, we will further to discuss the relevant conditions in setting both $f$ and $g$ are constant functions in subsection \ref{subsection2.4}.

Set  $d_{H}(T)=\sum_{x\in T}d_{H}(x)$ and $f(S)=\sum_{x\in S}f(x)$.
The lemma below will be used in the demonstration process of  our Theorem \ref{theorem1}-\ref{theorem3}.
\begin{lemma}\label{lemma1}{\rm (Gao \cite{Gao12})} Assume $G$ is a graph, functions $f$ and $g$ are
integer-valued on its vertex set meeting $f(x)\ge g(x)$ for every $x$ in $G$. Set $n'$,   $m\in\Bbb N^{+}\cup\{0\}$. Then $G$ is fractional $(g,f,n',m)$-critical deleted
iff
\begin{eqnarray}\label{1}
&\quad&f(S)-g(T)+d_{G-S}(T)\\\nonumber
&\ge&\max_{U\subseteq S,   H\subseteq
E(G-U), |U|=n', |H|=m}\{f(U)+\sum_{x\in T}d_{H}(x)-e_{H}(T,S)\}
\end{eqnarray}
for any subsets $S,T$ of $V(G)$ with $S\cap T=\emptyset$ and $|S|\ge n'$.
\end{lemma}

In very special circumstances, $n'$ vertices consist an independent set, then it comes to {\em fractional
ID-$(g,f,m)$-deleted graph}. Analogously, when functions $g$ and $f$ take special value for all vertices, it becomes different names which are presented in Table \ref{table3}.
\begin{table}
\caption{Special cases of fractional ID-$(g,f,m)$-deleted graphs}\label{table3}
\begin{tabular}{ll}
\hline
setting (for any $v\in V(G)$) & name \\
\hline
$g=f$ & fractional ID-$(f,m)$-deleted
graph\\
\hline
$f(x)=b$ and $g(x)=a$ &fractional
ID-$(a,b,m)$-deleted graph\\
\hline
$m=0$ &fractional
ID-$(g,f)$-factor-critical graph\\
\hline
\end{tabular}
\end{table}

The following results in fractional ID-$(g,f,m)$-deleted setting as second part of main conclusions
which are the extension of Theorem 4, Theorem 5 and Theorem 6
showed in Gao et al. \cite{Gaowei4}, respectively.
\begin{theorem}\label{theorem4}
Assume $G$ is a graph with $n$ vertices, and $b,a,m,\Delta$ are
non-negative integers meeting $b-\Delta\ge a\ge2$ and $n>
\frac{(b+\Delta+2a)(b+2m+a-2)}{\Delta+a}$. Functions $g,f$ are integer-valued
on its vertex set satisfy $b-\Delta\ge f(x)-\Delta\ge g(x)\ge a$ for every vertex $x$.  Then $G$ is fractional
ID-$(g,f,m)$-deleted if $\delta(G)\ge\frac{(b+a)n}{b+2a+\Delta}$.
\end{theorem}

\begin{theorem}\label{theorem5}
Assume $G$ is a graph with $n$ vertices, and $b,a,m,\Delta$ are
non-negative integers meeting $b-\Delta\ge a\ge2$,
$\delta(G)\ge\frac{(\Delta+a)n}{b+2a+\Delta}+\frac{b(b-\Delta)}{\Delta+a}+m$ and $n>
\frac{(b+a+2m-1)(b+2a+\Delta)}{\Delta+a}$. Functions $g,f$
as integer-valued on its vertex set satisfy $b-\Delta\ge f(x)-\Delta\ge g(x)\ge a$ for every vertex $x$ in $G$. Then $G$ is fractional ID-$(g,f,m)$-deleted if for any $xy\ne E(G)$, we have
$$\max\{d_{G}(y),d_{G}(x)\}\ge\frac{(b+a)n}{b+2a+\Delta}.$$
\end{theorem}

\begin{theorem}\label{theorem6}
Assume $G$ is a graph with $n$ vertices, and $b,\Delta,a,m$ as
non-negative integers meeting $b-\Delta\ge a\ge2$,
$\delta(G)\ge\frac{(\Delta+a)n}{b+2a+\Delta}+\frac{b(b-\Delta)}{a+\Delta}+m$ and $n>
\frac{(b+2m+a-2)(b+\Delta+2a)}{\Delta+a}$. Functions $g,f$
are integer-valued on its vertex set satisfy $b-\Delta\ge f(x)-\Delta\ge g(x)\ge a$ for every vertex $x$ in $G$. Then $G$ is
fractional ID-$(g,f,m)$-deleted if $\sigma_{2}(G)\ge\frac{2(b+a)n}{b+2a+\Delta}$.
\end{theorem}

\section{Proof of first part results: Theorem \ref{theorem1}-\ref{theorem3}}
By observing, we find that $\delta(G)\ge\frac{(b-\Delta)n+(a+\Delta)n'}{a+b}$ in
Theorem \ref{theorem1} implies
$\sigma_{2}(G)\ge\frac{2((b-\Delta)n+(a+\Delta)n')}{a+b}$ and
$\delta(G)\ge n'+m+\frac{(b-\Delta)b}{\Delta+a}$ in Theorem
\ref{theorem3}. Hence, it's sufficient to make Theorem
\ref{theorem2} and \ref{theorem3} proved.

We deduce the conclusion on graph without non-adjacent vertices below.

\begin{lemma}\label{lemma2} Assume $G$ is a complete graph with $n$ vertices, and
$b,\Delta,a, n',m$ as non-negative integers meeting $b-\Delta\ge a\ge 2$ and $n>\frac{(b+a+2m-2)(b+a)}{a+\Delta}+n'$.  Functions $g,f$ as integer-valued on its vertex set satisfy $b-\Delta\ge f(x)-\Delta\ge g(x)\ge a$ for every vertex $x$ in $G$. Then $G$ is
fractional $(g,f,n',m)$-critical deleted.
\end{lemma}
{\bf Proof.}
Assume $G$ meets the conditions of Lemma \ref{lemma2} without being fractional $(g,f,n',m)$-critical deleted. Clearly,
$T\ne\emptyset$. According to Lemma
\ref{lemma1} and the fact that $\sum_{x\in
T}d_{H}(x)-e_{H}(T,S)$ at most $2m$, subsets
$T$ and $S$ of $V(G)$ with $T\cap S\emptyset$ exist to satisfy
$$f(S)+d_{G-S}(T)-g(T)\le \max_{U\subseteq S,|U|=n'}f(U)-1+2m$$
or
\begin{equation}\label{2}
f(S-U)-g(T)+d_{G-S}(T)\le2m-1,
\end{equation}
in which $|S|\ge n'$. Choose $T$ and $S$ with minimum $|T|$.
Hence, for every $x\in T$, we derive $b-1-\Delta\ge g(x)-1\ge d_{G-S}(x)$.

Note that $G-S$ is also complete for each vertex subset $S$. Thus, for
disjoint subsets $T,S\subseteq V(G)$, we deduce
\begin{eqnarray*}&\quad& f(S-U)-g(T)-2m+d_{G-S}(T)\\
&\ge&(|S|-n')(\Delta+a)+\sum_{x\in T}d_{G-S}(x)-(b-\Delta)|T|-2m\\
&\ge&(|S|-n')(\Delta+a)-(n-|S|)(b-n-\Delta+1+|S|)-2m\\
&=&(b+a-2n+1)|S|+|S|^{2}-(b-\Delta)n+n^{2}-n-(a+\Delta)n'-2m.
\end{eqnarray*}
Regarding it as the function of $|S|$, we look into the following cases in view of the fact that $|S|$ is an integer.

{\bf Case 1.} $a+b$ is even. Since
$n>\frac{(b+a+2m-2)(b+a)}{\Delta+a}+n'$ and $b+a\ge4$, we have
\begin{eqnarray*}&\quad& (b+a+1-2n)|S|+|S|^{2}-(b-\Delta)n+n^{2}-2m-n-(a+\Delta)n'\\
&\ge&(n-\frac{a+b}{2})(b+a+1-2n)+(n-\frac{b+a}{2})^{2}-n-(b-\Delta)n+n^{2}-2m-(a+\Delta)n'\\
&=&(\Delta+a)n-\frac{b+a}{2}-(\frac{b+a}{2})^{2}-2m-(\Delta+a)n'\\
&>&(\frac{(b+a+2m-2)(b+a)}{\Delta+a}+n')(\Delta+a)-2m-\frac{b+a}{2}-(\frac{b+a}{2})^{2}-(\Delta+a)n'\\
&=&\frac{3}{4}(b+a)^{2}+(b+a-1)2m-\frac{5}{2}(b+a)\\
&\ge&\frac{3}{4}\cdot16-\frac{5}{2}\cdot4>0,\end{eqnarray*} which
contradicts (\ref{2}).

{\bf Case 2.} $b-a\equiv 1$ (mod 2).  By
$n>\frac{(b+a+2m-2)(b+a)}{\Delta+a}+n'$ and $b+a\ge5$, we get
\begin{eqnarray*}&\quad& (b+a+1-2n)|S|+|S|^{2}-(b-\Delta)n+n^{2}-n-(\Delta+a)n'-2m\\
&\ge& (n-\frac{a+b+1}{2})(b+a-2n+1)+(n-\frac{b+a+1}{2})^{2}\\
&\quad&-(b-\Delta)n+n^{2}-2m-n-(a+\Delta)n'\\
&=&(\Delta+a)n-2m-(\frac{b+a+1}{2})^{2}-(a+\Delta)n'\\
&>&(\Delta+a)(\frac{(b+a)(b+a+2m-2)}{\Delta+a}+n')-(\frac{b+a+1}{2})^{2}-(\Delta+a)n'-2m\\
&=&(b+a-1)2m+\frac{3}{4}(b+a)^{2}-\frac{1}{4}-\frac{5}{2}(b+a)\\
&\ge&\frac{3}{4}\cdot25-\frac{5}{2}\cdot5-\frac{1}{4}>0,\end{eqnarray*}
 a contradiction.

The proof of complete graph setting is done. \qed

Clearly, Lemma \ref{lemma2} is the extension of previous
conclusion on the complete graph marked in Lemma 2 of Gao et al. \cite{Gaowei4}.
By setting $n'=0$ in Lemma \ref{lemma2}, the
corollary present below will be employed in Section \ref{secondpartsection}.
\begin{corollary}\label{corollary4} Assume $G$ is a complete graph having $n$ vertices, and
$b,\Delta,a,m$ as non-negative integers meeting $b-\Delta\ge a\ge2$ and $n>\frac{(b+a)(b+a+2m-2)}{a+\Delta}$.  Functions $g,f$ are
integer-valued on its vertex set satisfy $b-\Delta\ge f(x)-\Delta\ge g(x)\ge a$ for every vertex $x$ in $G$. Therefore, $G$ is
fractional $(g,f,m)$-deleted.
\end{corollary}

Graph is supposed to be non-complete in what follows. From this
point of view, the degree condition
$\max\{d_{G}(y),d_{G}(x)\}\ge\frac{(b-\Delta)n+(a+\Delta)n'}{a+b}$ for every $xy\ne E(G)$ in Theorem
\ref{theorem2} and $\sigma_{2}(G)\ge\frac{2((b-\Delta)n+(a+\Delta)n')}{a+b}$ in
Theorem \ref{theorem3} are meaningful.

\subsection{Correctness of Theorem \ref{theorem2}}\label{subsection2.1}
Assume $G$ meets all the assumptions of Theorem \ref{theorem2} without being fractional $(g,f,n',m)$-critical deleted. It can be inferred
$|T|\ge1$, and there exist disjoint subsets $T,S\subseteq V(G)$ satisfies (\ref{2}) with $|S|\ge n'$. We have
$b-1-\Delta\ge g(x)-1\ge d_{G-S}(x)$ for all vertex $x$ in $T$ by means of
selecting $S$ and $T$ with minimum $|T|$.

Let $d_{1}=\min\{d_{G-S}(x):x\in T\}$. We deduce $b-1-\Delta\ge d_{1}\ge0$ and
$$f(S-U)+d_{G-S}(T)-g(T)\ge d_{1}|T|+(\Delta+a)(|S|-n')-(b-\Delta)|T|.$$
This implies
\begin{equation}\label{3}
2m-1\ge (|S|-n')(\Delta+a)-(b-\Delta-d_{1})|T|.
\end{equation}
We choose vertex $x_{1}$ in $T$  to meet $d_{G-S}(x_{1})=d_{1}$.

If $|T|\le b-\Delta$, in terms of (\ref{3}) and $|S|+d_{1}\ge
d_{G}(x_{1})\ge\delta(G)\ge\frac{b(b-\Delta)}{a+\Delta}+n'+m$, we
verify
\begin{eqnarray*}
&\quad&2m-1\\
&\ge& (|S|-n')(\Delta+a)+|T|(\Delta+d_{1}-b)\\
&\ge& (a+\Delta)(\frac{b(b-\Delta)}{a+\Delta}+n'+m-d_{1}-n')+(\Delta+d_{1}-b)(b-\Delta)\\
&=&(b-a-\Delta)d_{1}+am+\Delta(m-d_{1}-\Delta+b)\\
&\ge&2m,
\end{eqnarray*}
which gets contradicted. Thus, $|T|\ge
b+1-\Delta\ge 1+a$.

On the condition that
$T-N_{T}[x_{1}]\ne\emptyset$, set $d_{2}=\min\{d_{G-S}(x):x\in
T-N_{T}[x_{1}]\}$ and take vertex $x_{2}$ belongs to $T-N_{T}[x_{1}]$ such that
$d_{G-S}(x_{2})=d_{2}$. Hence, $d_{1}\le d_{2}\le b-\Delta-1$.
Since $b-\Delta-1\ge d_{G-S}(x)$ for any vertex $x$ in $T$ and $|T|\ge
b-\Delta+1$, $T-N_{T}[x_{1}]\ne\emptyset$, thus $x_{1}$, $x_{2}$
must be existed. Considering the non-adjacent vertices assumption,
we deduce
$$\frac{n'(a+\Delta)+n(b-\Delta)}{b+a}\le\max\{d_{G}(x_{1}),d_{G}(x_{2}))\}\le|S|+d_{2},$$
which reveals
\begin{equation}\label{4}
|S|\ge\frac{n'(a+\Delta)+n(b-\Delta)}{b+a}-d_{2}.
\end{equation}

In view of $b-\Delta-d_{2}>0$ and $n-|S|-|T|\ge0$, we
infer
\begin{eqnarray*}
&\quad&(n-|T|-|S|)(b-\Delta-d_{2})\\
&\ge&(\Delta+a)(|S|-n')+\sum_{x\in T}(d_{G-S}(x)-b+\Delta)+1-2m\\
&\ge&(d_{1}+\Delta-b)|N_{T}[x_{1}]|-2m+(\Delta+a)|S|+1\\
&\quad&+(d_{2}-b+\Delta)(|T|-|N_{T}[x_{1}]|)-(a+\Delta)n'\\
&=&(a+\Delta)|S|+(d_{1}-d_{2})|N_{T}[x_{1}]|+(d_{2}+\Delta-b)|T|-(a+\Delta)n'-2m+1\\
&\ge&(d_{1}-d_{2})(1+d_{1})+(a+\Delta)|S|+(\Delta+d_{2}-b)|T|-(a+\Delta)n'-2m+1.
\end{eqnarray*}

It follows that
\begin{equation}\label{5}
0\le
n(b-d_{2}-\Delta)-(b+a-d_{2})|S|+2m+(1+d_{1})(d_{2}-d_{1})+(a+\Delta)n'-1.
\end{equation}

Using (\ref{4}), (\ref{5}),  $n>
\frac{(a+b)(a+b+2m-1)}{a+\Delta}+n'$ and $d_{1}\le d_{2}\le b-1-\Delta$, we obtain
\begin{eqnarray*}
0&\le&(b-d_{2}-\Delta)n-(b+a-d_{2})(\frac{(a+\Delta)n'+(b-\Delta)n}{b+a}-d_{2})+(d_{1}+1)(d_{2}-d_{1})\\
&\quad&+(a+\Delta)n'-1+2m\\
&=&-nd_{2}\frac{a+\Delta}{b+a}+d_{2}\frac{(a+\Delta)n'}{a+b}+(b+a)d_{2}-d_{1}^{2}-d_{2}^{2}+d_{1}d_{2}+d_{2}-d_{1}+2m-1\\
&<&-d_{1}^{2}-d_{2}^{2}+d_{1}d_{2}+2d_{2}-d_{1}+2m(1-d_{2})-1.\\
\end{eqnarray*}

If $d_{2}=0$, then we have $d_{1}=d_{2}=0$. According to
(\ref{4}), we have $|S|\ge\frac{(b-\Delta)n+(a+\Delta)n'}{a+b}$ and $|T|\le n-|S|\le
 \frac{(a+\Delta)n-(a+\Delta)n'}{a+b}$. By $\sum_{x\in
T}d_{H}(x)-e_{G}(T,S)\le d_{G-S}(T)$, we yield
\begin{eqnarray*}&\quad&f(S-U)-g(T)+d_{G-S}(T)-(\sum_{x\in T}d_{H}(x)-e_{G}(T,S))\\
&\ge&(a+\Delta)\cdot(\frac{(b-\Delta)n+(\Delta+a)n'}{b+a}-n')-(b-\Delta)\cdot\frac{(\Delta+a)n-n'(\Delta+a)}{b+a}\\
&\quad&+e_{G}(T,S)+d_{G-S}(T)-\sum_{x\in T}d_{H}(x)\\
&\ge&0,
\end{eqnarray*}
a contradiction.

If $d_{2}\ge1$,  we infer
$$0<-d_{1}^{2}-d_{2}^{2}+d_{1}d_{2}+2d_{2}-d_{1}+2m(1-d_{2})-1\le-d_{2}^{2}+(d_{1}+2)d_{2}-d_{1}^{2}-d_{1}-1.$$
Let
$$h_{1}(d_{2})=-d_{2}^{2}+(d_{1}+2)d_{2}-d_{1}^{2}-d_{1}-1.$$
This implies,
$$\max\{h_{1}(d_{2})\}=h_{1}(\frac{d_{1}+2}{2})=-\frac{3}{4}d_{1}^{2}\le0,$$
which is a contradiction. Thus, we  complete the derivation for the correctness.\qed

\subsection{Correctness of Theorem \ref{theorem3}}\label{subsection2.2}
Assume $G$ meets all the assumptions of Theorem \ref{theorem3} without being fractional $(g,f,n',m)$-critical deleted.
Apparently, $|T|\ge1$ and there exist  $T,S\subseteq V(G)$ with $T\cap S=\emptyset$ satisfies (\ref{2}) with $|S|\ge n'$.
Selecting $T$ and $S$ with minimum $|T|$, we obtain $b-1-\Delta\ge g(x)-1-\Delta\ge d_{G-S}(x)$ for any vertex $x$ in $T$.

Set $d_{1}$, $d_{2}$, $x_{1}$ and $x_{2}$ as defined before.
Similarly as discussed in Section 2.1, we yield $d_{1}\le d_{2}\le
b-\Delta-1$, $|T|\ge b-\Delta+1$ and $x_{1}$, $x_{2}$ must be
existed.

By means of degree assumption, we arrive
$$\frac{2(n(b-\Delta)+n'(\Delta+a))}{b+a}\le\sigma_{2}(G)\le2|S|+d_{2}+d_{1},$$
which reveals
\begin{equation}\label{6}
|S|\ge\frac{(a+\Delta)n'+(b-\Delta)n}{a+b}-\frac{d_{2}+d_{1}}{2}.
\end{equation}

Using the consideration in Subsection \ref{subsection2.1}, (\ref{5}) holds as well.
In light of (\ref{5}), (\ref{6}),  $n>
\frac{(b+a-2+2m)(b+a)}{\Delta+a}+n'$ and $d_{1}\le d_{2}\le b-1-\Delta$, we derive
\begin{eqnarray*}
0&\le&(d_{1}+1)(d_{2}-d_{1})+n(b-d_{2}-\Delta)-(b+a-d_{2})(\frac{n'(\Delta+a)+n(b-\Delta)}{b+a}-\frac{d_{1}+d_{2}}{2})\\
&\quad&+(\Delta+a)n'-1+2m\\
&=&d_{2}\frac{(\Delta+a)n'}{b+a}-nd_{2}\frac{\Delta+a}{b+a}+(b+a)\frac{d_{1}+d_{2}}{2}-d_{1}^{2}-\frac{d_{2}^{2}}{2}+\frac{d_{1}d_{2}}{2}+d_{2}-d_{1}+2m-1\\
&<&-d_{2}(a+b-3)+\frac{a+b}{2}(d_{1}+d_{2})-d_{1}^{2}-\frac{d_{2}^{2}}{2}+\frac{d_{1}d_{2}}{2}-d_{1}+2m(1-d_{2})-1.\\
\end{eqnarray*}

The case for $d_{2}=0$ can be proved in the similar way as
Subsection \ref{subsection2.1}.

If $d_{2}\ge1$, then we verify
\begin{eqnarray*}
0&<&-d_{2}(a+b-3)+\frac{a+b}{2}(d_{1}+d_{2})-d_{1}^{2}-\frac{d_{2}^{2}}{2}+\frac{d_{1}d_{2}}{2}-d_{1}+2m(1-d_{2})-1\\
&\le&-\frac{d_{2}^{2}}{2}-d_{2}(\frac{a+b}{2}-3-\frac{d_{1}}{2})-d_{1}^{2}+(\frac{a+b}{2}-1)d_{1}-1.\\
\end{eqnarray*}
Let
$$h_{2}(d_{2})=-\frac{d_{2}^{2}}{2}-d_{2}(\frac{a+b}{2}-3-\frac{d_{1}}{2})-d_{1}^{2}+(\frac{a+b}{2}-1)d_{1}-1.$$
If $d_{2}$ can reach to $3+\frac{d_{1}}{2}-\frac{a+b}{2}$ (i.e.,
$3+\frac{d_{1}}{2}-\frac{a+b}{2}\ge1$), then
$$\max\{h_{2}(d_{2})\}=h_{2}(3+\frac{d_{1}}{2}-\frac{a+b}{2}),$$
and $d_{2}\le 1$ in terms of  $b\ge a\ge2$ and $d_{1}\le b-1$.
Hence, $(d_{1}, d_{2})=(0,1)$ or $d_{1}=d_{2}=1$. By $b\ge a\ge
2$, we get $h_{2}(d_{2})\le0$ for both
$(d_{1},d_{2})=(0,1)$ and $(d_{1},d_{2})=(1,1)$, a contradiction.

If $d_{2}$ can't take
$3+\frac{d_{1}}{2}-\frac{a+b}{2}-\frac{1}{a+b}$ as its value, then
we have
\begin{eqnarray*}
0&<&-\frac{d_{2}^{2}}{2}-d_{2}(\frac{a+b}{2}-3-\frac{d_{1}}{2})-d_{1}^{2}+(\frac{a+b}{2}-1)d_{1}-1\\
&\le&-\frac{d_{1}^{2}}{2}-d_{1}(\frac{a+b}{2}-3-\frac{d_{1}}{2})-d_{1}^{2}+(\frac{a+b}{2}-1)d_{1}-1\\
&=&-d_{1}^{2}+2d_{1}-1\le 0,
\end{eqnarray*}
which also gets contradicted. In result, Theorem \ref{theorem3}
is proven.\qed

 \subsection{Sharpness}\label{sharpness} In this part, we give an example to prove the sharpness of the degree conditions in Theorem \ref{theorem1}-\ref{theorem3} in some sense. That is to say, the minimal condition $\delta(G)\ge\frac{(b-\Delta)n+(a+\Delta)n'}{a+b}$ can't be changed to $\delta(G)\ge\frac{(b-\Delta)n+(a+\Delta)n'}{a+b}-1$; we can't replace $\max\{d_{G}(y),d_{G}(x)\}\ge\frac{n(b-\Delta)+n'(a+\Delta)}{a+b}$ by $\max\{d_{G}(y),d_{G}(x)\}\ge\frac{n(b-\Delta)+n'(a+\Delta)}{a+b}-1$ in Theorem \ref{theorem2}; and the degree sum condition $\sigma_{2}(G)\ge\frac{2((b-\Delta)n+(a+\Delta)n')}{a+b}$ in Theorem \ref{theorem3} can't be transferred to $\sigma_{2}(G)\ge\frac{2((b-\Delta)n+(a+\Delta)n')}{a+b}-1$.

Let $b=a+\Delta$, $G_{1}=K_{at+n'}$ be a complete graph, $G_{2}=(bt+1)K_{1}$, and $G=G_{1}\vee
G_{2}$, where $t\in\Bbb N$ is a large number which ensures the graph to meet
$\delta(G)\ge m+n'+\frac{b(b-\Delta)}{\Delta+a}$ and $n>
\frac{(b+a-2+2m)(b+a)}{\Delta+a}+n'$), so
$n=|G_{1}|+|G_{2}|=(a+b)t+1+n'$. Let
$a=g(x)$ and $b=a+\Delta=f(x)$ for every vertex $x$ in $G$. We have
$$\frac{n'(a+\Delta)+n(b-\Delta)}{b+a}
>\delta(G)=n'+at>\frac{n'(a+\Delta)+n(b-\Delta)}{b+a}-1,$$
$$\frac{n'(a+\Delta)+n(b-\Delta)}{b+a}
>\max\{d_{G}(y),d_{G}(x)\}
=n'+at>\frac{n'(a+\Delta)+n(b-\Delta)}{b+a}-1,$$
$$\frac{2(n'(a+\Delta)+n(b-\Delta))}{b+a}>\sigma_{2}(G)=2(at+n')
\ge\frac{2(n'(a+\Delta)+n(b-\Delta))}{b+a}-1.$$
 Let
 $T=V(G_{2})$ and $S=V(G_{1})$, we get
\begin{eqnarray*}
&\quad&f(S)-g(T)+d_{G-S}(T)-\max_{U\subseteq S, |U|=n', H\subseteq
E(G-U), |H|=m}\{f(U)-e_{H}(T,S)+\sum_{x\in
T}d_{H}(x)\}\\
&=&b|S|-a|T|-(a+\Delta)n'=-a<0.
\end{eqnarray*}
According to Lemma \ref{lemma1}, $G$ isn't
fractional ($g,f,n',m$)-critical deleted.

 \subsection{Specific case in setting $(g,f)=(a,b)$}\label{subsection2.4}
According to the techniques in the proof of Lemma
\ref{lemma2}, we infer a likely conclusion for a graph
without non-adjacent vertices.
\begin{lemma}\label{lemma3}
Assume $G$ is a complete graph having $n$ vertices, and $b,n',a,
m,\Delta$ are non-negative integers meeting $n>\frac{(b+a-2+2m)(b+a)}{\Delta+a}+n'$ where $b-\Delta\ge a\ge2$. Then $G$ is
fractional $(a,b,n',m)$-critical deleted.
\end{lemma}

We arrive the corollary below by setting $n'=0$ in Lemma \ref{lemma3}, which is a sufficient condition for a fractional $(a,b,m)$-deleted complete graph.
\begin{corollary}\label{corollary5}
Assume complete graph $G$  having $n$ vertices, and $b,a,m,\Delta$ are non-negative integers meeting
$n>\frac{(b+a-2+2m)(b+a)}{\Delta+a}$ where $b-\Delta\ge a\ge2$. Then $G$ is fractional
$(a,b,m)$-deleted.
\end{corollary}

Note that Lemma \ref{lemma3} and Corollary \ref{corollary5} here
are the extension results for the corresponding conclusions in Gao
et al. \cite{Gaowei4}.

Set $f(x)=b$, $g(x)=a$ for arbitrary vertex $x$ in $G$. The
 necessary and sufficient condition
is achieved from Lemma \ref{lemma1}.
 \begin{lemma}\label{lemma4} Assume $G$ is a graph, $b$, $a$, $n'$, and  $m$ are non-negative
integers meeting $b\ge a$.  Therefore, $G$ is fractional
$(a,b,n',m)$-critical deleted iff for any disjoint subsets $T,S\subseteq V(G)$ with $|S|\ge n'$, we have
\begin{equation}\label{7}
b|S|+d_{G-S}(T)-a|T|\ge\max_{ |H|=m}\{(a+\Delta)n'+\sum_{x\in
T}d_{H}(x)-e_{H}(T,S)\}.
\end{equation}
\end{lemma}

Using Lemma \ref{lemma3} and Lemma \ref{lemma4}, in view of the
approaches used in Subsection \ref{subsection2.1} and Subsection \ref{subsection2.2}, we deduce the degree conditions depicted in Table \ref{table4} in fractional $(a,b,n',m)$-critical
deleted setting, which are corresponding to Theorem \ref{theorem1}-\ref{theorem3}.
We omit the detailed proof.
\begin{table}
\caption{Degree conditions in fractional
$(a,b,n',m)$-critical deleted setting }\label{table4}
\begin{tabular}{lll}
\hline
order of graph & degree condition& additional condition \\
\hline
$n>
\frac{(b+a-2+2m)(b+a)}{\Delta+a}+n'$ & $\delta(G)\ge\frac{(a+\Delta)n'+(b-\Delta)n}{b+a}$& \\
\hline
 $n>
\frac{(b+a-2+2m)(b+a)}{\Delta+a}+n'$&$\max\{d_{G}(x),d_{G}(y)\}\ge\frac{(\Delta+a)n'+(b-\Delta)n}{b+a}$ & $\frac{b(b-\Delta)}{\Delta+a}+m+n'\le \delta(G)$\\
\hline
 $n>
\frac{(b+a-2+2m)(b+a)}{\Delta+a}+n'$& $\sigma_{2}(G)\ge\frac{2(n'(\Delta+a)+n(b-\Delta))}{b+a}$& $\frac{b(b-\Delta)}{\Delta+a}+m+n'\le \delta(G)$ \\
\hline
\end{tabular}
\end{table}

Again, three theorems above present the new extension versions of
Theorem 7-9 in Gao et al.
\cite{Gaowei4}, respectively. Moreover, the example  in Subsection \ref{sharpness} shows that these degree
conditions in Table \ref{table4} are tight.

In particular, by taking $n'=0$ in   Table \ref{table4}, the corresponding degree conditions in fractional
$(a,b,m)$-deleted setting are obtained in Table \ref{table5}.
\begin{table}
\caption{Degree conditions in fractional
$(a,b,m)$-deleted setting}\label{table5}
\begin{tabular}{lll}
\hline
order of graph & degree condition& additional condition\\
\hline
$n>
\frac{(b+a-2+2m)(b+a)}{\Delta+a}$ & $\delta(G)\ge\frac{(b-\Delta)n}{b+a}$& \\
\hline
$n>
\frac{(b+a-2+2m)(b+a)}{\Delta+a}$ & $\max\{d_{G}(x),d_{G}(y)\}\ge\frac{(b-\Delta)n}{b+a}$ & $\delta(G)\ge
\frac{b(b-\Delta)}{\Delta+a}+m$ \\
\hline
$n>
\frac{(b+a-2+2m)(b+a)}{\Delta+a}$ & $\sigma_{2}(G)\ge\frac{2n(b-\Delta)}{a+b}$ & $\delta(G)\ge
\frac{b(b-\Delta)}{\Delta+a}+m$\\
\hline
\end{tabular}
\end{table}

\section{Proof of second part results: Theorem \ref{theorem4}-\ref{theorem6}}\label{secondpartsection}
Since $\delta(G)\ge\frac{n(a+b)}{2a+\Delta+b}$ in Theorem \ref{theorem4}
implies $\delta(G)\ge\frac{(\Delta+a)n}{b+2a+\Delta}+m+\frac{(b-\Delta)b}{\Delta+a}$ and
$\sigma_{2}(G)\ge\frac{2n(b+a)}{2a+\Delta+b}$ in Theorem \ref{theorem6},
it is sufficient for the proof of Theorem \ref{theorem5}-\ref{theorem6}.

\subsection{Correctness of Theorem \ref{theorem5}-\ref{theorem6}}\label{subsection3.1}

Here, first let's prove Theorem \ref{theorem5}. Let $G'=G-I$ for arbitrary independent set
$I$.  The conclusion is deduced by making sure that
$G'$ meets Table \ref{table2} or Corollary
\ref{corollary4}.

If every two vertices has an edge in $G'$, we
obtain
$$|G'|\ge\frac{n(b+a)}{b+2a+\Delta}>\frac{(b+a-1+2m)(b+a)}{a+\Delta}
>\frac{(b+2m+a-2)(b+a)}{\Delta+a}.$$
The conclusion holds in light of Corollary \ref{corollary4}.

If $I$ only contain one vertex, we yield
$|V(G')|>\frac{(b+2m+a-1)(b+2a+\Delta)-\Delta-a}{\Delta+a}>\frac{(b+a-1+2m)(b+a)}{\Delta+a}$.
Hence, $\delta(G')\ge\frac{(b-\Delta)b}{a+\Delta}+m$ and
$$\max\{d_{G'}(u),d_{G'}(v)\}\ge\frac{|V(G')|(b-\Delta)}{b+a}=\frac{(n-1)(b-\Delta)}{b+a}$$
for any $uv\notin E(G')$. Hence, the result obtained in view of Table \ref{table2}.

If $|I|\ge2$ and $G'$ isn't complete. Applying
degree condition, we infer
$|V(G')|\ge\frac{(b+a)n}{2a+b+\Delta}>\frac{(b+a)(b+2m+a-1)}{\Delta+a}$. If
$\max\{d_{G'}(u),d_{G'}(v)\}<\frac{|V(G')|(b-\Delta)}{b+a}$ for some
$uv\notin E(G')$, we arrive
$$\frac{(|V(G')|+|I|)(b+a)}{b+2a+\Delta}\le
\max\{d_{G}(v),d_{G}(u)\}<\frac{|V(G')|(b-\Delta)}{b+a}+|I|,$$
which implies
$$|V(G')|<\frac{(a+\Delta)(b+a)}{a^{2}+\Delta(2a+\Delta)}|I|\le\frac{(b+a)(a+\Delta)}{a^{2}+\Delta(2a+\Delta)}\frac{n(\Delta+a)}{2a+\Delta+b}=\frac{(b+a)n}{2a+\Delta+b}.$$
It contradicts  $|I|\ge2$ and $\max\{d_{G}(v),d_{G}(u)\}\ge\frac{(b+a)n}{b+2a+\Delta}$. Thus,
$$\max\{d_{G'}(v),d_{G'}(u)\}\ge\frac{|V(G')|(b-\Delta)}{b+a}$$ for any
$uv\notin E(G')$. Further, we get $m+\frac{b(b-\Delta)}{\Delta+a}\le \delta(G')$ in view of $|I|\le\frac{(a+\Delta)n}{2a+b+\Delta}$ and
$\frac{n(\Delta+a)}{b+2a+\Delta}+m+\frac{b(b-\Delta)}{a+\Delta}\le \delta(G)$. Therefore,
the result is obtained from Table \ref{table2}.

Hence, we finish the proof of Theorem \ref{theorem5}. By means
of Table \ref{table2} and Corollary \ref{corollary4},
Theorem \ref{theorem6} can be checked in the similar techniques. We
omit the detailed procedure. \qed

\subsection{Tight of results}\label{tight}
To show the tight of Theorem \ref{theorem4}, Theorem \ref{theorem5} and Theorem
\ref{theorem6}, we need the following lemma follows from the corollary of
Lemma \ref{lemma1}.
\begin{lemma}\label{kkk} Assume $G$ is a graph, functions $g,f$ are integer-valued on its vertex set meeting $f(x)\ge g(x)$ for every vertex $x$ in $G$. Set $m\in\Bbb N^{+}\cup\{0\}$.
Then $G$ is fractional $(g,f,m)$-deleted iff for all disjoint subsets $T,S\subseteq V(G)$, we have
\begin{equation}
f(S)+d_{G-S}(T)-g(T)\ge\max_{
|H|=m}\{\sum_{x\in T}d_{H}(x)-e_{H}(T,S)\}.
\end{equation}
\end{lemma}

Let $b=a+\Delta$.
Take $G=(bt+1)K_{1}\vee K_{at}\vee(bt+1)K_{1}$,
where $t\in\Bbb N$ is a large number. Apparently,
$n=2+(2b+a)t$. Set $f(x)=b$ and $g(x)=a$ for any vertex $x$ in $G$. We have
$$\frac{(b+a)n}{b+\Delta+2a}>\delta(G)=(b+a)t+1>\frac{(b+a)n}{b+\Delta+2a}-1,$$
$$\frac{(b+a)n}{b+\Delta+2a}>\max\{d_{G}(u),d_{G}(v)\}=(b+a)t+1>\frac{(b+a)n}{b+\Delta+2a}-1,$$
$$\frac{2(b+a)n}{b+\Delta+2a}>\sigma_{2}(G)=2+2(b+a)t>\frac{2(b+a)n}{b+\Delta+2a}-1.$$
Let $I=(bt+1)K_{1}$. For $G'=K_{at}\vee(bt+1)K_{1}$, let
$S=K_{at}$ and $T=(bt+1)K_{1}$.
We confirm that $e_{H}(T,S)=\sum_{x\in T}d_{H}(x)$ for arbitrary $H\subseteq E(G')$ having $m$ edges. As a result,
$$f(S)+d_{G-S}(T)-g(T)-(\sum_{x\in T}d_{H}(x)-e_{H}(T,S))=b(at)-a(bt+1)=-a.$$
To sum up, $G$ isn't fractional ID-$(g,f,m)$-deleted due to Lemma \ref{kkk} and $G'$ isn't fractional $(g,f,m)$-deleted.

\subsection{Specific case in setting $(g,f)=(a,b)$}\label{lastsubsection}

The below degree conditions in Table \ref{table6} in setting $g(x)=a$ and $f(x)=b$ are derived in terms of Corollary
\ref{corollary5}, Table \ref{table5}, and the approaches
in Subsection \ref{subsection2.4} and Subsection \ref{subsection3.1}.
\begin{table}
\caption{Degree conditions in fractional
ID-$(a,b,m)$-deleted setting}\label{table6}
\begin{tabular}{lll}
\hline
order of graph & degree condition& additional condition\\
\hline
$n>
\frac{(b+2a+\Delta)(b+a-2+2m)}{\Delta+a}$ & $\delta(G)\ge\frac{n(b+a)}{2a+\Delta+b}$& \\
\hline
$n>
\frac{(b+2a+\Delta)(b+a+2m-1)}{\Delta+a}$ & $\max\{d_{G}(x),d_{G}(y)\}\ge\frac{(b+a)n}{b+2a+\Delta}$ & $\delta(G)\ge\frac{(\Delta+a)n}{b+2a+\Delta}+\frac{b(b-\Delta)}{\Delta+a}+m$\\
\hline
$n>
\frac{(b+a-2+2m)(b+2a+\Delta)}{\Delta+a}$ & $\sigma_{2}(G)\ge\frac{2(b+a)n}{b+2a+\Delta}$ & $\delta(G)\ge\frac{(\Delta+a)n}{b+2a+\Delta}+\frac{b(b-\Delta)}{a+\Delta}+m$\\
\hline
\end{tabular}
\end{table}

One important thing we emphasize here is that the results presented in Table \ref{table6} are also the extensions of Theorem 10-12 in Gao et. al. \cite{Gaowei4}. Moreover, in terms
of the example presented in Subsection \ref{tight}, we
ensure that these degree conditions in Table \ref{table6} are also tight.

\section{Conclusion}
In our work, we mainly discuss the degree conditions for the existence of fractional factor in the setting that $b-\Delta\ge f(x)-\Delta\ge g(x)\ge a$ for each vertex $x$ in $G$, and some elements of graph are forbidden. Our results reveal that $\Delta$ is a key factor in this setting, and it specifically points out how $\Delta$ plays a role in the conclusion.

\section{Acknowledgments}
This research is partially supported by NSFC (Nos. 11761083, 11771402, 11671053).

\end{document}